\documentclass[a4paper,11pt]{article}
\usepackage{amsmath,amsthm,amssymb,amscd,epsfig}

\newtheorem{theo}{Theorem}[section]

\newtheorem{lem}[theo]{Lemma}

\newtheorem{coro}[theo]{Corollary}

\newtheorem{quest}[theo]{Question}
\newtheorem{conj}[theo]{Conjecture}

\newtheorem{rem}[theo]{Remark}

\numberwithin{equation}{section}

\def\R{\mathbb R}
\def\C{\mathbb C}
\def\N{\mathbb N}
\def\D{\mathbb D}
\def\Z{\mathbb Z}

\def\H{\mathcal H}
\def\M{\mathcal M}
\def\F{\mathcal F}

\def\diam{\text{diam}}
\def\Lip{\text{Lip}}

\title{Distortion of Hausdorff measures and improved Painlev\'e removability for  quasiregular mappings}

\author{\textit{K. Astala}\and \textit{A. Clop}\and \textit{J. Mateu} \and \textit{J. Orobitg} \and \textit{I. Uriarte-Tuero}
\thanks{{Astala was supported in part by the Academy of
Finland, projects 34082 and 41933. Clop, Mateu, Orobitg were supported by projects MTM2004-00519, HF2004-0208, 2005-SGR-00774. Uriarte-Tuero was supported by the Finnish Academy of Sciences}\newline
\newline Keywords   Quasiconformal, Hausdorff measure, Removability}}

\begin{document}

\maketitle

\begin{abstract}
The classical Painlev\'e theorem tells that sets of zero length are removable for bounded analytic functions, while (some) sets of  positive length are not.  For  general  $K$-quasiregular mappings in planar domains the corresponding critical dimension is $\frac{2}{K+1}$. We show that when $K>1$, unexpectedly  one has  improved removability. More precisely, we prove that sets $E$ of  $\sigma$-finite Hausdorff $\frac{2}{K+1}$-measure are removable  for bounded $K$-quasiregular mappings. On the other hand, $\dim(E) = \frac{2}{K+1}$ is not enough to guarantee this property.

We also study absolute continuity properties of pull-backs of Hausdorff measures under $K$-quasiconformal mappings, in particular at the relevant dimensions $1$ and $\frac{2}{K+1}$. For general Hausdorff measures ${\cal H}^t$,  $0 < t < 2$ ,  we reduce the absolute continuity properties  to an open question on conformal mappings, see Conjecture \ref{question1}.
\end{abstract}
\section{Introduction}

A homeomorphism $\phi:\Omega\rightarrow\Omega'$ between planar domains $\Omega,\Omega' \subset \C$ is called {\it{$K$-quasiconformal}} if  it belongs to  the Sobolev space $W^{1,2}_{loc}(\Omega)$ and satisfies the {\it{distortion inequality}}
\begin{equation}\label{distortioninequality}
\max_\alpha|\partial_\alpha\phi|\leq K\min_\alpha|\partial_\alpha\phi| \,\,\,\,\,\text{a.e. in }\Omega
\end{equation}
It has been known since the work of Ahlfors \cite{Al} that quasiconformal mappings preserve sets of zero Lebesgue measure. It is also well known that they preserve sets of zero Hausdorff dimension, since $K$-quasiconformal mappings are H\"older continuous with exponent $1/K$, see \cite{Mo}. However, these maps do not preserve Hausdorff dimension in general, and  it was in the work of the first author \cite{A} where the precise bounds for the distortion of   dimension were given. For any compact set $E$ with dimension $t$ and for any $K$-quasiconformal mapping $\phi$ we have
\begin{equation}\label{distortionofdimension}
\frac{1}{K}\left(\frac{1}{t}-\frac{1}{2}\right)\leq\frac{1}{\dim(\phi(E))}-\frac{1}{2}\leq K\left(\frac{1}{t}-\frac{1}{2}\right)
\end{equation}
Furthermore, these bounds are optimal, that is, equality may occur in either estimate. 

The fundamental question we study in this work is whether  the estimates (\ref{distortionofdimension})  can be improved to the level of Hausdorff measures $\H^t$. In other words, if $\phi$ is a planar $K$-quasiconformal mapping, $0<t<2$ and $t'=\frac{2Kt}{2+(K-1)t}$, we ask whether  it is  true that
\begin{equation}\label{abscont}
\H^t(E)=0\,\,\,\Rightarrow\,\,\,\H^{t'}(\phi(E))=0,
\end{equation}
or put briefly, $\phi^\ast\H^{t'}\ll\H^t$. Note that the above classical results of Ahlfors and Mori assert that this is true when $t=0$ or $t=2$. In fact \cite{A}, for the Lebesgue measure one has  even precise  quantitative bounds  $$|\phi(E)|\leq C\,|E|^\frac{1}{K}, $$
a result which also leads to the sharp Sobolev regularity,  $\phi\in W^{1,p}_{loc}(\C)$ for every  $p<\frac{2K}{K-1}$.

As a first main result of this paper  we prove (\ref{abscont}) for $t=\frac{2}{K+1}$, i.e. for the case of image dimension  $t'=1$. 

\begin{theo}\label{firstmaintheo}
Let $\phi$ be a planar $K$-quasiconformal mapping, and let $E$ be a compact set. Then,
\begin{equation}\label{festimate}
\M^1(\phi(E))\leq C\,\left(\M^\frac{2}{K+1}(E)\right)^\frac{K+1}{2K}
\end{equation}
As a consequence, 
$$\H^\frac{2}{K+1}(E)=0\,\,\,\Rightarrow\,\,\,\H^{1}(\phi(E))=0$$
\end{theo}

Here $\M^t$ denotes $t$-dimensional Hausdorff content. As one of the key points in proving Theorem \ref{firstmaintheo} we show that for  planar quasiconformal mappings $h$ which are {\it{conformal}} in the complement  $\C\setminus E$  the inequality (\ref{festimate}) improves strongly: such mappings $h$ essentially preserve  the
$1$-dimensional Hausdorff content of the compact set $E$,
\begin{equation}\label{confoutside}
\M^1(h(E))\leq C_K\,\M^1(E)
\end{equation}
The constant  $C_K$ depends only on $K$ if $h$ is normalized at $\infty$, requiring  $h(z) = z + {\cal O}(1/z)$.
For the  area  the corresponding estimate  was shown in \cite{A}.  In fact, as we will see later, a counterpart of (\ref{confoutside}) for the $t$-dimensional Hausdorff content $\M^t$ is the only missing detail for proving the absolute continuity  $\phi^\ast\H^{t'}\ll\H^t$ for general $t$. 
Towards solving (\ref{abscont}) we conjecture  that actually
$$\M^t(h(E))\leq C\,\M^t(E), \quad \quad 0 < t \leq 2
$$
whenever $E \subset \C$ compact and $h$ is normalized and conformal in $\C \setminus E$  admitting  a $K$-quasiconformal extension to $\C$. For a more detailed  discussion and other formulations see  Section 2. 

The reason our methods work only in the special case of dimension $t=1$  is that the  content $\M^1$ is equivalent to a suitable $BMO$-capacity \cite{V}. For dimensions  $1<t<2$, we do have interpolating estimates but unfortunately we have to settle for  the Riesz capacities. We have
$$
C_{\alpha, t}(h(E))\leq C\,C_{\alpha, t}(E)
$$
for any $t\in (1,2)$ and  $\alpha=\frac{2}{t}-1$. This fact has  consequences towards the absolute continuity of Hausdorff measures under quasiconformal mappings, but these bounds are not strong enough for (\ref{abscont}) when $1<t'<2$.\\

Recall that $f$ is a {\it{$K$-quasiregular mapping}} in a domain $\Omega\subset\C$ if  $f\in W^{1,2}_{loc}(\Omega)$ and $f$ satisfies the distortion inequality (\ref{distortioninequality}).
 When $K=1$, this class agrees with the class of analytic functions on $\Omega$.
The classical {\it{Painlev\'e problem}} consists of giving metric and geometric characterizations of those sets $E$ that are removable for bounded analytic functions. Here Painlev\'e's theorem tells us that sets of zero length are removable, while Ahlfors \cite{Ahl2} showed that no set of Hausdorff dimension  $> \;1$ has this property. For the related $BMO$-problem  Kaufman  \cite{K} proved that the condition $\H^1(E)=0$ is a precise characterization for removable singularities of  $BMO$ analytic functions. Thus for analytic removability, dimension $1$ is the critical point both for  $L^\infty$ and  $BMO$. However, the solution to the original Painlev\'e problem lies much deeper and was only recently achieved by Tolsa (\cite{T1},\cite{T2}) in terms of curvatures of measures. Under the assumption that $\H^1(E)$ is finite, Painlev\'e's problem was earlier solved by G. David \cite{D}, who showed that a set $E$ of positive and finite length is removable for bounded analytic functions if and only if it is purely unrectifiable. Furthermore, the countable semiadditivity of analytic capacity, due to Tolsa \cite{T1}, asserts that this result remains true if we only assume $\H^1(E)$ to be $\sigma$-finite.\\
\\
It is now natural to approach the Painlev\'e problem for $K$-quasiregular mappings. We say that a compact set $E$ is {\it{removable for bounded $K$-quasiregular mappings}}, or simply {\it{$K$-removable}}, if for every open set $\Omega\supset E$, every bounded $K$-quasiregular mapping $f:\Omega\setminus E\rightarrow\C$ admits a $K$-quasiregular extension to $\Omega$. In this definition, as in the analytic setting, we may replace $L^\infty(\Omega)$ by $BMO(\Omega)$ to get a close variant of the problem.\\
\\
The sharpness of the bounds in equation (\ref{distortionofdimension}) determines the index $\frac{2}{K+1}$ as the critical dimension in both the $L^\infty$ and $BMO$ quasiregular removability problems. In fact, Iwaniec and Martin previously conjectured \cite{IM}  that in $\R^n$, $n \geq 2$, sets with Hausdorff measure  $\H^\frac{n}{K+1}(E)=0$ are removable for bounded $K$-quasiregular mappings. A preliminary positive answer for $n=2$ was described  in   \cite{AIM}. Generalizing this, in the present work we show that  surprisingly, for $K>1$ one can do even better: 
we have the following improved Painlev\'e removability.

\begin{theo}\label{improvePainleve}
Let $K > 1$ and suppose $E$ is any compact set with $$\H^\frac{2}{K+1}(E) \quad \; \sigma-\mbox{finite}. $$
Then $E$ is removable for all bounded $K$-quasiregular mappings. 
\end{theo}
\noindent The theorem fails for $K = 1$, since for instance the line segment $E = [0,1]$  is not removable. 

For  the converse direction, the work \cite{A} finds for every  $t> \frac{2}{K+1}$ non-$K$-removable sets with $\dim(E) = t$. We make an improvement also here and construct compact sets with dimension precisely equal to $\frac{2}{K+1}$ yet not removable for some  bounded $K$-quasiregular mappings. For details see Theorem \ref{construction}. \\

The  above theorems  are  closely connected via  the classical  Stoilow factorization, which tells \cite{AIM},  \cite{LV} that in planar domains $K$-quasiregular mappings are precisely the maps $f$ representable in the form $f = h \circ \phi$,  where $h$ is analytic and $\phi$ is $K$-quasiconformal.
Indeed, the first step in proving  Theorem \ref{improvePainleve}  will be to show that for a general $K$-quasiconformal mapping $\phi$ one has
$$\H^\frac{2}{K+1}(E)\,\,\,\,\text{$\sigma$-finite}\,\,\,\Rightarrow\,\,\,\H^{1}(\phi(E))\,\,\,\,\text{$\sigma$-finite}$$
However, this conclusion will not be enough since there are rectifiable sets of finite length, such as  $E = [0,1]$, that are non-removable for bounded analytic functions. Therefore, in addition,  we  need to establish that such  'good'  sets of  positive analytic capacity 
actually  behave better also under quasiconformal mappings. That is, we show that up to a set of zero length, 
$$F  \quad 1\mbox{-rectifiable} \,\,\,\Rightarrow\,\,\,\dim(\phi(F))>\frac{2}{K+1}$$
For details and a precise formulation see Corollary \ref{maintheorem4}.\\
\\
The paper is structured as follows. In Section 2 we deal with the quasiconformal distortion of Hausdorff measures and of other set functions. 
In Section 3 we study the quasiconformal distortion of $1$-rectifiable sets. Section 4 gives the proof for the improved Painlev\'e removability theorem for $K$-quasiregular mappings and other related questions. Finally in   section 5 we describe a construction of non-removable sets.

\section{Absolute Continuity}
 
There are several natural ways to normalize the quasiconformal mappings $\phi:\C \to\C$.
In this work we mostly  use the {\it{principal}} $K$-quasiconformal mappings, i.e. mappings that are conformal outside a compact set  and are normalized by $\phi(z)-z=O\left(\frac{1}{|z|}\right)$ as $|z|\to\infty$.

It is shown in the work \cite{A} of the first author that  for all  $K$-quasiconformal mappings $\phi:\C\rightarrow\C$, 
\begin{equation}\label{above}
|\phi(E)|\leq C\,|E|^{1/K} 
\end{equation}
where $C$ is a constant that depends on the normalizations. By scaling we may always arrange 
\begin{equation}\label{norma}
\diam(\phi(E)) = \diam(E) \leq 1
\end{equation}
and then $C=C(K)$ depends only on $K$. In order to achieve the result (\ref{above}),  one first reduces to the case where the set $E$ is a finite union of disks. Secondly, applying Stoilow factorization methods the mapping $\phi$ is written as $\phi=h\circ \phi_1$, where both $h,\phi_1:\C\rightarrow\C$ are $K$-quasiconformal mappings, such that $\phi_1$ is conformal on $E$ and $h$ is conformal in the complement of the  set  $F=\phi_1(E)$. Here one obtains the right conclusion for $\phi_1$, 
$$|\phi_1(E)|\leq C\,|E|^\frac{1}{K}$$
by including $\phi_1$ in a holomorphic family of quasiconformal mappings. Further, one shows in \cite[p. 50]{A} that under the special assumption where $h$  is conformal outside of $F$, we have
\begin{equation}\label{areathing}
|h(F)|\leq C\,|F|
\end{equation}
where the constant $C$ still depends only on $K$. 

In searching for absolute continuity properties of other Hausdorff measures under quasiconformal mappings, such a decomposition seems unavoidable, and this leads one to look for counterparts of (\ref{areathing}) for Hausdorff measures $\H^t$ or Hausdorff contents $\M^t$. Here we establish the result for the dimension $t=1$.

\begin{lem}\label{gamma0invariance}
Suppose $E\subset\C$ is a compact set, and let $\phi:\C\rightarrow\C$ be a principal $K$-quasiconformal mapping, such that $\phi$ is conformal on $\C\setminus E$. Then,
$$\M^1(\phi(E))\simeq \M^1(E)$$
with constants depending only on $K$.
\end{lem}

In order to prove this result  some background is needed. The space of functions of bounded mean oscillation, $BMO$, is invariant under quasiconformal changes of variables \cite{R}. More precisely, if $\phi$ is a $K$-quasiconformal mapping and $f\in BMO(\C)$, then $f\circ\phi\in BMO(\C)$ with $BMO$-norm
$$\|f\circ\phi\|_\ast\leq C(K)\,\|f\|_\ast$$
The space $BMO(\C)$ gives rise to a capacity,
$$\gamma_0(F)=\sup |f'(\infty)|$$
where the supremum runs over all functions $f\in BMO(\C)$ with $\|f\|_ \ast\leq 1$, that are holomorphic on
$\C\setminus E$ and satisfy  $f(\infty)=0$. Here $f'(\infty) = \lim_{| z|\to \infty} z \left(f(z)-f(\infty)\right)$.
Observe that in this situation $\overline\partial f$ defines a  distribution supported on $F$, and actually
$|\langle\overline\partial f, 1\rangle|=|f'(\infty)|$. It turns out  \cite{V} that for any compact set $E$ we have
\begin{equation}\label{kaufman}
\gamma_0(E)\simeq\M^1(E).
\end{equation}
According to the theorem of Kaufman   \cite{K}, in the class of  functions  $f\in BMO(\C)$  holomorphic on
$\C\setminus E$ every  $f$ admits a holomorphic extension to the  whole plane if and only if $\M^1(E)=0$.
That is, $\gamma_0$ characterizes those compact sets which are  removable for $BMO$ holomorphic functions.
Because of these equivalences,  to prove Lemma \ref{gamma0invariance}  it suffices to show that $\gamma_0(\phi(E))\simeq \gamma_0(E)$. 

\begin{proof}[Proof of Lemma \ref{gamma0invariance}]
Suppose that $f\in BMO(\C)$ is a holomorphic mapping of $\C\setminus E$ such that $\|f\|_\ast\leq 1$ and $f(\infty)=0$. Then the function $g=f\circ\phi^{-1}$ is  in $BMO(\C)$ and $\|g\|_{\ast}\leq C(K)$. On the other hand, $g$ is holomorphic on $\C\setminus \phi(E)$, and since $\phi$ is a principal $K$-quasiconformal mapping, $g(\infty)=0$ and 
$$|g'(\infty)|=\lim_{|z|\to\infty}|z\,g(z)|=\lim_{|w|\to\infty}|\phi(w)\,f(w)|=|f'(\infty)|.$$
Hence, $\gamma_0(E)\leq C(K)\,\gamma_0(\phi(E))$. The converse inequality follows  by symmetry, since also the inverse $\phi^{-1}$ is a principal mapping,.
\end{proof}

This lemma is a first step towards the results  on absolute continuity, as presented in the following reformulation of Theorem \ref{firstmaintheo}.

\begin{theo}\label{zerotozero}
Let $E$ be a compact set and $\phi:\C \to \C$ $K$-quasiconformal, normalized by (\ref{norma}). Then
$$\M^1(\phi(E))\leq C\,\left(\M^\frac{2}{K+1}(E)\right)^\frac{K+1}{2K}$$
where the constant $C =C(K)$  depends only on $K$. In particular, if $\H^\frac{2}{K+1}(E)=0$ then $\H^1(\phi(E))=0$.
\end{theo}
\begin{proof} There is no restriction if we assume $E\subset\D$. We can also assume that $\phi$ is a principal $K$-quasiconformal mapping, conformal outside  $\D$. Now, given $\varepsilon>0$, there is a finite covering of  $E$ by open disks $D_j=D(z_j,r_j)$, $j=1,...,n$, such that
$$\sum_{j=1}^nr_j^\frac{2}{K+1}\leq\M^\frac{2}{K+1}(E)+\varepsilon$$
Denote  $\Omega=\cup_{j=1}^nD_j$. As in \cite{A}, we use a decomposition $\phi=h\circ\phi_1$, where both $\phi_1,h$ are principal $K$-quasiconformal mappings. 
Moreover, we may require that $\phi_1$ is conformal in $\Omega\cup(\C\setminus\D)$ and that $h$ is conformal outside $\phi_1(\overline{\Omega})$.\\
\\By Lemma \ref{gamma0invariance}, we see that
$$\M^1(\phi(E))\leq\M^1(\phi(\Omega))=\M^1(h\circ\phi_1(\Omega))\leq C\,\M^1(\phi_1(\Omega))$$
Hence the problem has been reduced to estimating  $\M^1(\phi_1(\Omega))$. For this, $K$-quasidisks have area comparable to the square of the diameter,
$$\diam(\phi_1(D_j))\simeq|\phi_1(D_j)|^{1/2}=\left(\int_{D_j}J(z,\phi_1)\,dA(z)\right)^\frac{1}{2}$$
with constants which depend only on $K$. Thus, using H\"older estimates twice, we obtain 
$$\sum_{j=1}^n\diam(\phi_1(D_j))\leq C(K)\,\left(\sum_{j=1}^n \int_{D_j}J(z,\phi_1)^pdA(z)\right)^\frac{1}{2p}\left(\sum_{j=1}^n|D_j|^\frac{p-1}{2p-1}\right)^{1-\frac{1}{2p}}$$
as long as $J(z,\phi_1)^p$ is integrable.
But here we are in the special situation of \cite[Lemma 5.2]{AN}. Namely, as $\phi_1$ is conformal in the subset $\Omega$, we may  take $p=\frac{K}{K-1}$ and apply \cite{AN} to obtain
$$
\sum_{j=1}^n \int_{D_j}J(z,\phi_1)^pdA(z) = \int_{\Omega}J(z,\phi_1)^pdA(z)\leq \pi
$$
With the above choice of $p$ one has $\frac{p-1}{2p-1}=\frac{1}{K+1}$. Hence we get
\begin{equation}\label{conformalabscont}
\sum_{j=1}^n\diam(\phi_1(D_j))\leq C(K)\,\left(\sum_{j=1}^nr_j^\frac{2}{K+1}\right)^\frac{K+1}{2K}\leq C(K)\left(\M^\frac{2}{K+1}(E)+\varepsilon\right)^\frac{K+1}{2K}
\end{equation}
But $\cup_j\phi_1(D_j)$ is a covering of $\phi_1(\Omega)$, so that actually we have
$$\M^1(\phi(E))\leq C \M^1(\phi_1(\Omega))\leq C(K)\left(\M^\frac{2}{K+1}(E)+\varepsilon\right)^\frac{K+1}{2K}$$
Since this holds for every $\varepsilon>0$, the result follows.
\end{proof}

At this point we want to emphasize that for a general quasiconformal mapping $\phi$ we have $J(z,\phi)\in L^p_{loc}$ only for $p<\frac{K}{K-1}$. The improved borderline integrability ($p=\frac{K}{K-1}$) under the extra assumption that $\phi_{\big | \Omega}$ is conformal  was shown in \cite[Lemma 5.2]{AN}. This phenomenon was crucial for our argument, since we are studying Hausdorff measures rather than dimension. Actually, the same procedure shows that inequality (\ref{conformalabscont}) works in a much more general setting.  That is,  still under the special assumption that $\phi_1$ is conformal in $\cup_{j=1}^nD_j$, we have for any $t\in [0,2]$
\begin{equation}\label{abscontconfasumption}
\left(\sum_{j=1}^n\diam(\phi_1(D_j))^{d}\right)^\frac{1}{d}\leq C(K)\,\left(\sum_{j=1}^n\diam(D_j)^t\right)^{\frac{1}{t}\,\frac{1}{K}}
\end{equation}
where $d=\frac{2Kt}{2+(K-1)t}$. On the other hand, another  key point in our proof was the estimate 
$$\M^1(h(E))\leq C\,\M^1(E),$$
valid whenever $h$ is a principal $K$-quasiconformal mapping which is conformal outside  $E$. We believe that finding the counterpart to this estimate is crucial for understanding distortion of Hausdorff measures under quasiconformal mappings. We make the following

\begin{conj}\label{question1}
Suppose we are given a real number $d\in (0,2]$. Then for any  compact set $E \subset \C$ and for any principal $K$-quasiconformal mapping $h$ which is conformal on $\C\setminus E$, we have
\begin{equation}\label{question}
\M^d(h(E))\simeq\M^d(E)
\end{equation}
with constants that  depend on $K$ and $d$ only.
\end{conj}

\noindent One may also formulate a convenient discrete variant, which  is actually stronger than Conjecture \ref{question1}.

\begin{quest} \label{kysymys}
Suppose we are given a real number $d\in (0,2]$ and a finite number of disjoint disks $D_1,...,D_n$. If a  mapping $h$ is conformal on $\C\setminus \cup_{j=1}^nD_j$ and admits a  $K$-quasiconformal extension to $\C$, is it then true that
\begin{equation} 
\sum_{j=1}^n\diam\left(h(D_j)\right)^d \simeq  \sum_{j=1}^n\diam(D_j)^d
\end{equation}
with constants that depend only on $K$ and $d$ ?
\end{quest}

We already know that (\ref{question}) is true for $d=1$ and $d=2$; however for Question \ref{kysymys} we know a proof only  at $d=2$. An affirmative answer to Conjecture \ref{question1}, combined  with the optimal integrability bound proving (\ref{abscontconfasumption}), would provide the absolute continuity of $\phi^\ast\H^d$ with respect to $\H^t$, where $d=\frac{2Kt}{2+(K-1)t}$,  $0\leq t\leq 2$ and $K\geq 1$. Therefore, (\ref{question}) would have important consequences in the theory of quasiconformal mappings. \\
\\
The positive answer to (\ref{question}) for the dimension $d=1$  was based on the equivalence (\ref{kaufman}) and  the invariance of $BMO$. Actually more is true:  the space $VMO$,  equal to the $BMO$-closure of uniformly continuous functions, is  quasiconformally invariant as well. We may also describe  $VMO$  as consisting  of  functions $f\in BMO$ for which
$$\lim\frac{1}{|B|}\int_B|f-f_B|=0$$
as $|B|+\frac{1}{|B|}\to\infty$. As we now  see, the invariance of VMO has interesting consequences. 

\begin{theo}\label{finitetofinite}
Let $E\subset\C$ be a compact set, and $\phi:\C\rightarrow\C$ a $K$-quasiconformal mapping. If $\H^\frac{2}{K+1}(E)$ is finite (or even $\sigma$-finite), then $\H^1(\phi(E))$ is $\sigma$-finite.
\end{theo}

This result may be equivalently expressed in terms of the lower Hausdorff content. To understand this alternative formulation of Theorem \ref{finitetofinite}, we first need some background. A {\it{measure function}} is a continuous non-decreasing function $h(t)$, $t\geq 0$, such that $\lim_{t\to 0}h(t)=0$. If $h$ is a measure function and $F\subset\C$ we set
$$\M^h(F)=\inf \sum_jh(\delta_j)$$
where the infimum is taken over all countable coverings of $F$ by disks of diameter $\delta_j$. When $h(t)=t^\alpha$, $\alpha>0$, $\M^h(F)=\M^\alpha(F)$ equals  the $\alpha$-dimensional Hausdorff content of $F$. Moreover, the content $\M^\alpha$ and the measure $\H^\alpha$ have the same zero sets.
We will denote by $\F =\F_d$ the class of measure functions $h(t)=t^d\,\varepsilon(t)$, $0\leq\varepsilon(t)\leq 1$, such that $\lim_{t\to 0}\varepsilon(t)=0$. The lower $d$-dimensional Hausdorff content of $F$ is then defined by
$$\M^d_{\ast}(F)=\sup_{h\in\F_d} \M^h(F)$$
One has $\M^d_\ast\leq \M^d$ but it can happen that $\M^d_\ast(F)=0<\M^d(F)$. For instance, if $F$ is the segment $[0,1]$ in the plane, then $\M^1_\ast(F)=0$ but $\M^1(F)=1$. An old result of Sion and Sjerve \cite{SS} in geometric measure theory asserts that $\M^d_\ast(F)=0$ if and only if $F$ is a countable union of sets with finite $d$-dimensional Hausdorff measure. For a disk $B$, $\M^d_\ast(B)=\M^d(B)$, and for open sets $U$, $\M^d_\ast(U)\simeq \M^d(U)$. We may now reformulate Theorem \ref{finitetofinite} as follows.

\begin{theo}\label{finitetofinitealternative}
Let $E\subset\C$ be a compact set, and $\phi:\C\rightarrow\C$ a principal $K$-quasiconformal mapping. If $\M^\frac{2}{K+1}_\ast(E)=0$, then $\M^1_\ast(\phi(E))=0$. 
\end{theo}

For the proof, for any bounded set $F\subset\C$ define first
\begin{equation}\label{gammah}
\gamma_{\ast}(F)=\sup|f'(\infty)|
\end{equation}
where the supremum is taken over all functions $f\in VMO$, with $\|f\|_\ast\leq 1$, which are holomorphic on $\C\setminus F$ and satisfy  $f(\infty)=0$. Again here we may replace $|f'(\infty)|$ by $|\langle\overline\partial f, 1\rangle|$. The $VMO$ invariance leads to the following analogue of Lemma 
\ref{gamma0invariance}.

\begin{lem}\label{gammavmoinvariance}
Let $E$ be a compact set. For any principal $K$-quasiconformal mapping $\phi:\C\rightarrow\C$, conformal on $\C\setminus E$, we have
$$\gamma_\ast(\phi(E))\simeq\gamma_{\ast}(E).$$
\end{lem}
\begin{proof}
Consider $f\in VMO$ which is analytic in $\C\setminus\phi(E)$ and $f(\infty)=0$. Set $g=f\circ\phi$. Then $g\in VMO$, $g$ is analytic on $\C\setminus E$, $\|g\|_\ast\leq C\,\|f\|_\ast$ and $|g'(\infty)|=|f'(\infty)|$ since $\phi$ is a principal $K$-quasiconformal mapping. Consequently $\gamma_\ast(\phi(E))\leq C\,\gamma_{\ast}(E)$.
\end{proof}

It was shown by Verdera  that this $VMO$ capacity is essentially the $1$-dimensional lower content.

\begin{lem}[{\cite{V}, p. 288}]\label{gammavmolowercontent}
For any compact set $E$, $\M^1_\ast(E)\simeq\gamma_{\ast}(E)$.
\end{lem}

With these tools we are ready to prove Theorem \ref{finitetofinitealternative}. 

\begin{proof}[Proof of Theorem \ref{finitetofinitealternative}]
Naturally, the argument is similar to that in Theorem \ref{zerotozero}. Without loss of generality, we may assume that $E\subset\D$ and that $\phi$ is a principal $K$-quasiconformal mapping. Furthermore, we may assume that  $\H^\frac{2}{K+1}(E)$ is finite, and for any $\delta$ we have  a finite family of disks $D_i$ such that $E\subset\cup_iD_i$, $\sum_i\diam(D_i)^\frac{2}{K+1}\leq\H^\frac{2}{K+1}(E)+1$ and $\diam(D_i)<\delta$. Set $\Omega=\cup_iD_i$. Again, we have a decomposition $\phi=\phi_2\circ\phi_1$, where both $\phi_1$ and $\phi_2$ are principal $K$-quasiconformal mappings, and where we may require that $\phi_1$ is conformal in $(\C\setminus\overline\D)\cup\Omega$, and $\phi_2$ is conformal outside $\phi_1(\overline\Omega)$. Thus,
$$\M^1_\ast(\phi(E))\leq\M^1_\ast(\phi(\overline\Omega))$$
By Lemma \ref{gammavmolowercontent}, the lower content can be replaced by the $VMO$ capacity,
$$\M^1_\ast(\phi(\overline\Omega))\leq C\,\gamma_\ast(\phi(\overline\Omega))$$
Since $\phi_2$ is conformal outside of $\phi_1(\overline{\Omega})$, from Lemma \ref{gammavmoinvariance} we obtain
$$\gamma_\ast(\phi(\overline\Omega)) = \gamma_\ast(\phi_2\circ\phi_1(\overline\Omega))\simeq\gamma_\ast(\phi_1(\overline\Omega))  \leq C\,\M^1_\ast(\phi_1(\overline\Omega))$$
where the last inequality uses again  Lemma \ref{gammavmolowercontent}. It hence remains to estimate  $\M^1_\ast(\phi_1(\overline{\Omega}))$. For this,  take $h\in\F$, $h(t)=t\,\varepsilon(t)$, and argue as in Theorem \ref{zerotozero}. Since $K$-quasiconformal mappings are H\"older continuous with exponent $1/K$,
$$
\aligned
\M^h(\phi_1(\overline\Omega))
&\leq \; \sum_j\diam(\phi_1(D_j))\varepsilon(\diam(\phi_1(D_j)))
\;  \leq \; \varepsilon(C_K\,\delta^{1/K})\,\sum_j\diam(\phi_1(D_j))\\
&\leq\varepsilon(C_K\,\delta^{1/K})\,\sum_j\left(\int_{D_j}J(z,\phi_1)^\frac{K}{K-1}\,dm(z)\right)^\frac{K-1}{2K}|D_j|^{1/2K}\\
&\leq\varepsilon(C_K\,\delta^{1/K})\,C_K\,\left(\sum_j\diam(D_j)^\frac{2}{K+1}\right)^\frac{K+1}{2K}
\leq\varepsilon(C_K\,\delta^{1/K})\,C_K\,\left(\H^\frac{2}{K+1}(E)+1\right)^\frac{K+1}{2K}
\endaligned
$$
Finally, taking $\delta\to 0$ we get $\M^h(\phi(E))=0$. This holds for any $h\in\F$, and the Theorem follows.
\end{proof}

One might think  of extending  the preceeding results from the critical index $\frac{2}{K+1}$ to arbitrary ones by using other capacities that behave like a Hausdorff content. For instance, the capacity $\gamma_\alpha$, associated to analytic functions with the $\Lip(\alpha)$ norm \cite{OF}, satisfies
$$\M^{1+\alpha}(E)\simeq\gamma_{\alpha}(E)$$
but unfortunately, the space $\Lip(\alpha)$ is not invariant under a quasiconformal change of variables. Thus, other procedures are needed. It turns out that the homogeneous Sobolev spaces provide suitable tools, basically since $\dot{W}^{1,2}(\C)$ is  invariant under quasiconformal mappings. Here recall that for $0<\alpha<2$ and $p>1$, the homogeneous Sobolev space $\dot{W}^{\alpha,p}(\C)$ is defined as the space of Riesz potentials 
$$f=I_\alpha\ast g$$
where $g\in L^p(\C)$ and $I_\alpha(z)=\frac{1}{|z|^{2-\alpha}}$. The norm is given by $\|f\|_{\dot{W}^{\alpha,p}(\C)}=\|g\|_p$. When $\alpha=1$,  $\dot{W}^{1,p}(\C)$ agrees with the space of functions $f$ whose first order distributional derivatives are given by $L^p(\C)$ functions. Let $f\in\dot{W}^{1,2}(\C)$ and let $\phi$ be a $K$-quasiconformal mapping on $\C$. Defining  $g=f\circ\phi$ we have
$$
\aligned
\int_\C|Dg(z)|^2dA(z)
&\leq \int_\C|Df(\phi(z))|^2\,|D\phi(z)|^2\,dA(z)\\
&\leq\,K\,\int_\C|Df(\phi(z))|^2\,J(z,\phi)\,dA(z)\\
&=K\int_\C|Df(w)|^2\,dA(w)
\endaligned
$$
so that $g\in\dot{W}^{1,2}(\C)$. In other words, every $K$-quasiconformal mapping $\phi$ induces a bounded linear operator
$$T:\dot{W}^{1,2}(\C)\rightarrow\dot{W}^{1,2}(\C), \quad \quad T(f)=f\circ\phi$$
with norm depending only on $K$. As we have mentioned before, this operator $T$ is also bounded on $BMO(\C)$ \cite{R}. Moreover, Reimann and Rychener \cite[p.103]{RR} proved that $\dot{W}^{\frac{2}{q},q}(\C)$, $q>2$, may be represented as a complex interpolation space between $BMO(\C)$ and $\dot{W}^{1,2}(\C)$. It follows that $T$ is bounded on the Sobolev spaces $\dot{W}^{\frac{2}{q},q}(\C), q>2$. More precisely, there exists a constant $C=C(K,q)$ such that
\begin{equation}\label{rychener}
\|f\circ\phi\|_{\dot{W}^{\frac{2}{q},q}(\C)}\leq C \|f\|_{\dot{W}^{\frac{2}{q},q}(\C)}
\end{equation}
for any $K$-quasiconformal mapping $\phi$ on $\C$. These invariant function spaces provide us with related invariant capacities.
Recall (e.g. \cite[pp.34 and 46]{AH}) that for any pair $\alpha>0$, $p>1$ with $0<\alpha p< 2$, one defines the Riesz capacity of a compact set $E$ by
$$\dot{C}_{\alpha, p}(E)=\sup\{\mu(E)^p\}$$
where the supremum runs over all positive measures $\mu$ supported on $E$, such that $\|I_\alpha\ast\mu\|_q\leq 1$, $\frac{1}{p}+\frac{1}{q}=1$. We get an equivalent capacity if we replace positive measures $\mu$ by distributions $T$ supported on $E$, $\|I_\alpha\ast T\|_q\leq 1$, and take the supremum of $|\langle T,1\rangle|^p$.

To see the connection with equation (\ref{rychener})  consider the set functions
$$\gamma_{1-\alpha,q}(E)=\sup\{|f'(\infty)|; f \mbox{ analytic in } \C\setminus E, \|f\|_{\dot{W}^{1-\alpha,q}}\leq 1 \mbox{  and } f(\infty)=0\}$$
Observe again that $|f'(\infty)|=|\langle\overline\partial f, 1\rangle|$ where this action must be understood in the sense of distributions. With this terminology we have

\begin{lem}
Suppose that $E$ is a compact subset of the plane. Then, for any $p\in(1,2)$,
$$\dot{C}_{\alpha, p}(E)^{1/p}\simeq\gamma_{1-\alpha,q}(E)$$
where $\alpha=\frac{2}{p}-1$ and $q=\frac{p}{p-1}$.
\end{lem}
\begin{proof}
On one hand, let $\mu$ be an admissible measure for $\dot{C}_{\alpha, p}$. Then, $I_\alpha\ast\mu$ is in $L^q$ with norm at most $1$. Define $f=\frac{1}{z}\ast\mu$. Clearly, $f$ is analytic outside $E$, $f(\infty)=0$ and $f'(\infty)=\mu(E)$. Moreover, up to multiplicative constants, 
$$\widehat{f}(\xi)\simeq\frac{1}{\xi}\,\widehat{\mu}(\xi)=\frac{\overline\xi}{|\xi|}\,\frac{1}{|\xi|}\,\widehat{\mu}(\xi)=\widehat{R}\,\widehat{I_1}\,\widehat{\mu}$$
and consequently we can write
$$f=\frac{1}{z}\ast\mu=R(I_1\ast\mu)=I_{1-\alpha}\ast R(I_\alpha\ast\mu)$$
where $R$ is a Calder\'on-Zygmund operator and $\|f\|_{\dot{W}^{1-\alpha,q}}=\|R(I_\alpha\ast\mu)\|_q\lesssim\|I_\alpha\ast\mu\|_q$.

For the converse, let $f=I_{1-\alpha}\ast g$ be an admissible function for $\gamma_{1-\alpha, q}$. We have that, up to a multiplicative constant, 
 $T=\overline\partial f$ is  an admissible distribution for $\dot{C}_{\alpha, p}$ because 
$$I_\alpha\ast T=R^t (g)$$
where $R^t$ is the transpose of $R$.
Thus, $\dot{C}_{\alpha, p}(E)^{1/p}\geq |\langle T,1\rangle|=|f'(\infty)|$ and the proof is complete.
\end{proof}

We end up with new quasiconformal invariants built on the Riesz capacities.

\begin{theo}\label{invariantrieszcapacity}
Let $\phi:\C\rightarrow\C$ be a principal $K$-quasiconformal mapping of the plane, which is conformal on $\C\setminus E$. Let $1< p< 2$ and $\alpha=\frac{2}{p}-1$. Then
$$\dot{C}_{\alpha,p}(\phi(E))\simeq\dot{C}_{\alpha,p}(E)$$
with constants that depend only on $K$ and $p$.
\end{theo}
\begin{proof} By the preceding Lemma it  suffices to show that $\gamma_{1-\alpha, q}\left(\phi (E)\right) \leq C_K \gamma_{1-\alpha, q}(E)$.

Let $f$ be an admissible function for $\gamma_{1-\alpha, q}(\phi(E))$. This means that $f$ is holomorphic on $\C\setminus\phi(E)$, $f(\infty)=0$ and that $\|f\|_{\dot{W}^{1-\alpha,q}}\leq 1$. Then, we consider the function $g=f\circ\phi$. Clearly, $\overline\partial (f\circ\phi) = 0$  outside $E$ and $g(\infty)=0$. Moreover, for $\alpha=\frac{2}{p}-1$ we have $1-\alpha=\frac{2}{q}$. Hence, because of equation (\ref{rychener}), 
$$\|g\|_{\dot{W}^{1-\alpha,q}}\leq C_K \|f\|_{\dot{W}^{\frac{2}{q},q}}\leq C(K,q)$$
so that $\frac{1}{C(K,q)}\,g$ is an admissible function for $\gamma_{1-\alpha, q}(E)$. Hence, as $\phi$ is a principal $K$-quasiconformal mapping,
$$\gamma_{1-\alpha,q}(E)\geq\frac{1}{C(K,q)} |f'(\infty)|$$
and  we may take supremum over $f$. \end{proof}

The above theorem  has direct  consequences towards the absolute continuity of Hausdorff measures, but unfortunately these are slighly weaker than one would wish for. In fact, there are compact sets $F$ such that $C_{\alpha, p}(F)=0$ and $\H^h(F)>0$, for some measure function $h(t)=t^p\,\varepsilon(t)$. Thus, Theorem \ref{invariantrieszcapacity} does not help for Conjecture \ref{question1}.
We have to content with the following setup:

Given $1<d<2$ consider the measure functions $h(t)=t^d\,\varepsilon(t)$ where
\begin{equation} \label{morestuff}
\int_0 \varepsilon(t)^\frac{1}{d-1}\frac{dt}{t}<\infty
\end{equation}
Typical examples of such functions are $h(t) = t^d |\log t|^{-s}$ or $h(t) = t^d |\log t|^{1-d} \log^{-s}(|\log t|)$ where $s> d-1$.

\begin{coro}
Let $E$ be a compact set on the plane, and $\phi:\C\rightarrow\C$ a principal $K$-quasiconformal mapping, conformal outside of $E$. Let $1<d<2$. Then,
$$\M^h(\phi(E))\leq C\,\M^d(E)$$
for any measure function $h(t)=t^d\,\varepsilon(t)$ satisfying (\ref{morestuff}). 
Moreover, if $\H^d(E)<\infty$ then $\H^h(\phi(E))=0$ for every such $h$.
\end{coro}
\begin{proof}
By \cite[Theorem 5.1.13]{AH}, given a measure function $h$ satisfying   (\ref{morestuff}) there is a constant  $C=C(h)$ with
$$\M^h(\phi(E))\leq C\,C_{\alpha, d}(\phi(E)), \hskip20pt \alpha=\frac{2}{d}-1$$
By Theorem \ref{invariantrieszcapacity},  $C_{\alpha, d}(\phi(E))\leq C\,C_{\alpha, d}(E)$ and  using again  \cite[Theorem 5.1.9]{AH} we  finally have
$C_{\alpha,d}(E)\leq \M^d(E)$.
\end{proof}

Arguing now as in Theorems \ref{zerotozero} and \ref{finitetofinitealternative}, we arrive at the following conclusion. 
\begin{coro}
Let $E$ be a compact set of the plane and suppose $\phi:\C\rightarrow\C$ is a $K$-quasiconformal mapping. Let $t\in(\frac{2}{K+1},2)$ and $d=\frac{2Kt}{2+(K-1)t}$. Then, under the normalization 
(\ref{norma}),
$$
\M^h(\phi(E)) \leq C\; C_{\alpha, d}(\phi(E)) \leq C\,\left(\M^t(E)\right)^\frac{1}{Kt}, \quad \alpha = \frac{2}{d}-1
$$
for any measure function $h$ satisfying (\ref{morestuff}). The constant $C$ depends only on $h$ and 
$K$.
\end{coro}

\noindent Here note that  for $\frac{2}{K+1}<t<2$ we always have $1<d<2$ in the above Corollary.

\section{Distortion of rectifiable sets}

In general, if $\phi$ is a $K$-quasiconformal mapping and $E$ is a compact set, it follows from (\ref{distortionofdimension}) that
\begin{equation}\label{kasi}
\dim(E) = 1 \quad  \Rightarrow \quad \frac{2}{K+1}\leq\dim\phi(E)\leq\frac{2K}{K+1}
\end{equation}
Here for both estimates one may find mappings $\phi$ and sets $E$ such that the 
equality is  attained, see \cite{A}.
There all  examples   come from  non regular Cantor-type constructions. Thus the extremal distortion of Hausdorff dimension is attained, at least, by sets irregular enough. The main purpose of this section is to prove that some irregularity is also necessary. Namely, we show that  quasiconformal images of $1$-rectifiable sets cannot achieve the maximal distortion of dimension.

\begin{theo}\label{maintheorem}
Suppose that $\phi:\C\rightarrow\C$ is a $K$-quasiconformal mapping. Let $E\subset\partial\D$ be a subset of the unit circle with  $\dim(E)=1$. Then we have the strict inequality
$$\dim(\phi(E))>\frac{2}{K+1}$$
\end{theo}
\smallskip

\noindent With similar but easier argument one may also prove that for such sets $E$, neither can $\dim(\phi(E))$ attain the upper bound in (\ref{kasi}). For details see Remark \ref{yla}.

From this Theorem we obtain as an immediate corollary the following more general result.

\begin{coro}\label{maintheorem4}
Suppose that $E$ is a $1$-rectifiable set, and let $\phi:\C\rightarrow\C$ be a $K$-quasiconformal mapping. Then there exists a subset $E_0\subset E$ of zero length such that 
$$\dim\phi(E\setminus E_0)>\frac{2}{K+1}.$$
\end{coro}

\noindent  Recall that a set $E\subset\C$ is said to be {\it{$1$-rectifiable}} if there exists a set $E_0$ of zero length such that $E\setminus E_0$ is contained in a countable union of Lipschitz curves, that is,
$$E\setminus E_0\subset\bigcup_{j=1}^\infty\Phi_j([0,1])$$ 
where all $\Phi_j:[0,1]\rightarrow\C$ are Lipschitz  mappings. Alternatively \cite{M} $1$-rectifiable sets can be viewed as subsets countable unions of ${\cal C}^1$ curves, modulo a set of zero length. In particular, for any $\varepsilon>0$ there is a decomposition
$$E\setminus E_0'=\bigcup_{i=1}^\infty E_i$$
where $E_0'$ has zero length and each $E_i$ can be written as $E_i=f_i(F_i)$, with $f_i:\C\rightarrow\C$ a $(1+\varepsilon)$-bilipschitz mapping and $F_i\subset\partial\D$. From this and Theorem \ref{maintheorem} we obtain Corollary \ref{maintheorem4}. \\
\\
To prove Theorem \ref{maintheorem}, first some reductions may be made. Recall \cite{LV} that every $K$-quasiconformal mapping $\phi$ can be factored as $\phi=\phi_n\circ\dots\circ\phi_1$ where each $\phi_j$ is $K_j$-quasiconformal, and $K_1\,K_2\dots\cdot K_n=K$. In particular, given $\varepsilon>0$, we can choose $K_j\leq 1+\varepsilon$ for all $j=1,...,n$, when $n$ is large enough. On the other hand, recall that from the distortion of Hausdorff dimension (\ref{distortionofdimension}) we have
\begin{equation}\label{minimaldimension}
\frac{1}{\dim\phi(E)}-\frac{1}{2}\leq K\left(\frac{1}{\dim E}-\frac{1}{2}\right)
\end{equation}
If $\phi$ is such that equality in (\ref{minimaldimension}) holds for $E$, then every factor $\phi_j$ above must  give equality for the set $E_j = \phi_{j-1}\circ\dots\phi_1(E)$ and $K=K_j$. In other words, if one of the mappings $\phi_j$ fails to satisfy the equality in (\ref{minimaldimension}), then so will $\phi$.  By combining these facts, we deduce that in order to prove Theorem \ref{maintheorem} we can assume that $K=1+\varepsilon$ with $\varepsilon>0$ as small as we wish. \\
\\

For  mappings with small dilatation it is possible achieve 
quantitative and more symmetric local distortion estimates.  In particular,  Theorem \ref{maintheorem} will follow from the next lower bounds for compression of dimension.

\begin{theo} \label{ensi}
Suppose $\phi:\C \to \C$ is $(1+\varepsilon)$-quasiconformal and $E\subset \partial \D$. Then for all $\varepsilon>0$ small enough,
\begin{equation}\label{squareddistortion}
\dim(E)\geq 1-c_0\,\varepsilon^2\,\,\,\Rightarrow\,\,\,\dim(\phi(E))\geq 1-c_1 \,\varepsilon^2 
\end{equation}
where the  constants  $c_0, \, c_1 >0$ are independent of  $\varepsilon$.
\end{theo}

 Our basic strategy towards this result  is to reduce it  to the properties of harmonic measure and 
 conformal  mappings admitting quasiconformal extensions. Indeed, denote by $\mu$ the Beltrami coefficient of  $\phi$ and let $h$ be  the principal solution to $\overline\partial h=\chi_{\D}\mu\,\partial h$.
Then $h$ is conformal outside the unit disk. Inside   $\D$ it has the same dilatation $\mu$ as $\phi$, and hence differs from this by a conformal factor. Consequently, we may find  Riemann mappings $f: \D \to \Omega : = \phi(\D)$  and $g: \D \to \Omega' : = h(\D)$  so that 
\begin{equation} \label{factors}
 \phi(z) = f \circ g^{-1} \circ h (z), \quad z \in \D
\end{equation}
Moreover, since the $(1+\varepsilon)$-quasiconformal mapping $G = g^{-1} \circ h$ preserves the disk, reflecting across the boundary  $\partial \D$ one may extend  $G$  to a $(1+\varepsilon)$-quasiconformal mapping of $\C$. At the same time,  this procedure provides both $f$ and $g$  with   $(1+\varepsilon)^2$-quasiconformal extensions to the entire  plane $\C$.

As the final reduction we  now find from (\ref{factors})  that  for Theorems \ref{maintheorem} and \ref{ensi}  it is sufficient  to prove the following result.
 
\begin{theo}\label{mainlemma}
Suppose that $f:\C\rightarrow\C$ is a $({1+\varepsilon})$-quasiconformal mapping of $\C$, conformal in the disk $\D$. Let $A\subset\partial\D$. 
There are constants $c_0$, $c_1$ and $\gamma_0$, $\gamma_1$, independent of $\varepsilon$,
such that for $\varepsilon \geq 0$ small enough,
$$(i) \quad \quad  \dim(A) \geq 1-c_0\,\varepsilon^2\,\,\,\Rightarrow\,\,\,\dim(f(A)) \geq 1-c_1\,\varepsilon^2$$
and
$$(ii) \quad \quad \dim(A) \leq 1-\gamma_0\,\varepsilon^2\,\,\,\Rightarrow\,\,\,\dim(f(A)) \leq 1-\gamma_1\,\varepsilon^2$$
\end{theo}

\begin{proof}
The first conclusion $(i)$ follows from  Makarov's fundamental estimates for the harmonic measure \cite{Ma}, see also  \cite[p.231]{P}. In the work  \cite{Ma} Makarov proves  that for any conformal mapping $f$ defined on $\D$, for any Borel subset $A\subset\partial\D\,$ and for every $q>0$ we have the lower bound
\begin{equation}\label{makarov}
\dim(f(A))\geq\frac{q\,\dim(A)}{\beta_f(-q)+q+1-\dim(A)}
\end{equation} 
Here  $\beta_f(p)$ stands for  the integral means spectrum. That is, for a given $p\in \R$, $\beta_f(p)$ is  the infimum of all numbers $\beta$ such that
\begin{equation}\label{beta}
\int_0^{2\pi}|f'(re^{it})|^pdt=O\left(\frac{1}{(1-r)^\beta}\right)
\end{equation}
as $r\to 1^{-}$. 

We hence need  estimates for $\beta_f(p)$, and here for mappings admitting  $K$-quasiconformal extensions one has qualitively sharp bounds. Indeed, it can be shown \cite[p.182]{P} that
\begin{equation}\label{pommerenke}
\beta_f(p)\leq 9\left(\frac{K-1}{K+1}\right)^2\,p^2
\end{equation}
for any $p\in\R$.
The constant $9$ is not optimal but suffices for our purposes. Choosing $q=1$ in (\ref{makarov}) gives immediately the first claim $(i)$. \\
\\
For general conformal mappings there is no  bound for expansion of dimension, i.e. there is no  upper bound analogue of  (\ref{makarov}). Hence the proof of $(ii)$ uses strongly the fact that mappings considered have $({1+\varepsilon})$-quasiconformal extensions. However, also here this information is easiest to use in the form (\ref{pommerenke}).

 We first need to introduce some further notation. The Carleson squares of the unit disk are defined as
$$Q_{j,k}=\big\{z\in\D:2^{-k}\leq 1-|z|< 2^{-k+1}, \;2^{-k+1}\pi j\leq\arg(z)< 2^{-k+1}\pi(j+1)\big\}$$
Given a point $z\in\D\setminus\{0\}$, let $Q(z)$  denote the unique Carleson square that contains $z$. 
Then it  follows  from Koebe's distortion Theorem and quasisymmetry  \cite{AIM}, \cite{LV}  that if $D(\xi,r)$ is a disk centered at $\xi\in\partial D$, we have
\begin{equation}\label{koebe2}
\diam(f (D))\simeq \diam(f (Q(z)))\simeq |f'(z)|(1-|z|),  \quad \quad {\rm for} \;\; z=(1-r)\xi,
\end{equation}
whenever $f:\C\rightarrow\C$ is a $K$-quasiconformal mapping, conformal in $\D$.\\
\\
Furthermore, assume  we are given a family of  disjoint disks $D_i=D(\xi_i, r_i)$ with centers $\xi_i\in\partial\D$, $i\in\N$, on the unit circle. Write then $z_i=(1-r_i)\xi_i$, and for any pair of real numbers $0<\alpha<\delta<1$  define two subsets of indices,
$$\aligned
I_g(\alpha,\delta)&=\big\{i\in\N; |f'(z_i)|\leq(1-|z_i|)^{\frac{\alpha}{\delta}-1}\big\}\\
I_b(\alpha,\delta)&=\N\setminus I_g(\alpha,\delta)
\endaligned$$\\
Diameter sums over the 'good' indexes $I_g(\alpha,\delta)$ are easy to estimate. We have
$$
\sum_{i\in I_g(\alpha, \delta)}\diam(f(D_i))^\delta\leq\,C\,\sum_{i\in I_g(\alpha, \delta)}|f'(z_i)|^\delta\,(1-|z_i|)^\delta\leq\,C\,\sum_{i\in I_g(\alpha, \delta)}(1-|z_i|)^\alpha
$$
where $C$ depends only on $K$. In other words, 
\begin{equation}\label{football}
\sum_{i\in I_g(\alpha,\delta)}\operatorname{diam}(f(D_i))^\delta\leq C\,\sum_{i\in I_g(\alpha,\delta)}\operatorname{diam}(D_i)^\alpha
\end{equation}
\\
It is well known that the integral means can be used to control   the complementary indexes $I_b(\alpha,\delta)$. We give the technical details in a separate Lemma:

\begin{lem}\label{hardpart}
Assume that $0 < \alpha = 1-M\varepsilon^2$, for some $M>400$, and let $\delta = \alpha( 1+N\varepsilon^2) $, where $20\sqrt{M} < N < M$. Then 
$$\sum_{i\in I_b(\alpha,\delta)}\operatorname{diam}(f(D_i))^\delta\leq C$$
where $C$ is independent of $D_j$. Moreover, $\delta$ satisfies   $\delta < 1 - \gamma \varepsilon^2$ where
$\gamma = M-N > 0$.
\end{lem}
\begin{proof}
We classify the bad indexes $I_b(\alpha,\delta)$ by defining for $k=1,2,...$, and $m\in\Z$
$$
I^k_m=\big\{i\in I_b(\alpha,\delta); \; 2^{-k}\leq 1-|z_i|<2^{1-k}, \; \;2^{-1-m}\leq|f'(z_i)|(1-|z_i|)\leq 2^{-m}\big\}
$$
and write $q^k_m=\#I^k_m$. By (\ref{koebe2})  $|f'(z_i)|\,(1-|z_i|) $  is comparable to $\diam(f(D_i))$, which is always smaller than $\diam(f(3\D))$. On the other hand, if $i\in I^k_m$ then
$$(2^{-k})^\frac{\alpha}{\delta}\leq(1-|z_i|)^\frac{\alpha}{\delta}<(1-|z_i|)|f'(z_i)|\leq 2^{-m}$$
Hence the indexes $m$ with $I^k_m$ nonempty  lie on an interval $m_0 \leq 
m\leq\frac{\alpha}{\delta}k$. 

From Koebe we also see that  if $i\in I^k_m$ then $|f'(w)|^p \sim 2^{p(k-m)}$ for every  $w\in Q(z_i)$, with constants depending only on $p$. Combining this with (\ref{pommerenke}) gives for any $\tau > 0$
$$q^k_m\leq C\,2^{k\left(\frac{9}{4}\varepsilon^2p^2+\tau+1-\frac{k-m}{k}p\right)}$$
where $C$ now depends on $p$ and $\tau$. We may take $p=\frac{k-m}{10k\varepsilon^2}$ and obtain
$$q^k_m\leq\,C\,2^{(1+\tau)k-\frac{k}{(10\varepsilon)^2}\left(\frac{k-m}{k}\right)^2}$$
Since  $\diam(f(D_i))$ is comparable to $|f'(z_i)|\,(1-|z_i|) \sim 2^{-m}$ for $i\in I_m^k$,
\begin{equation}\label{sum}
\sum_{i\in I_b(\alpha,\delta)}\diam(f(D_i))^\delta\leq C \sum_{k=0}^\infty\sum_{m = m_0}^{\frac{\alpha}{\delta}k} q_m^k2^{-m\delta}\leq C \sum_{k=0}^\infty\sum_{m= m_0}^{\frac{\alpha}{\delta}k} 2^{k\left(1+\tau-m\frac{\delta}{k}-\frac{1}{100\varepsilon^2}\left(\frac{k-m}{k}\right)^2\right)}\\
\end{equation}
One  now needs to ensure that the exponent $1+\tau-m\frac{\delta}{k}-\frac{1}{100\varepsilon^2}\left(\frac{k-m}{k}\right)^2$  is negative. In particular, we want the exponent to attain its maximum at $m=\frac{\alpha}{\delta}k$, and this is satisfied if 
$$\frac{\alpha}{\delta}\leq 1-\frac{1}{2}(10\varepsilon)^2\,\delta$$
Under the assumptions of the Lemma this is easy to verify. Similarly one verifies that  the specific choices of the Lemma yield the  maximum value 
$$
1+\tau-\alpha-\frac{1}{(10\varepsilon)^2}\left(1-\frac{\alpha}{\delta}\right)^2 < 0
$$
 when  $\tau$ is chosen small enough. It follows that the sum in (\ref{sum}) has a finite upper bound depending only on the   constants $M$, $N$.
 This proves Lemma \ref{hardpart}.
\end{proof}

The  dimension bounds required in  part $(ii)$ of  Theorem \ref{mainlemma} are now easy to establish.
For every $\alpha > 1-\gamma_0 \varepsilon^2$ we have  coverings of $A$, consisting of   families  of disks $D_j=D(z_j, r_j)$ centered on  $\partial\D$ and radius $r_j \leq \rho \to 0$ uniformly small, such that the sums
$\sum_j\diam(D_j)^\alpha$
are uniformly bounded. On the image side, for each $\delta>0$
$$\sum_i\diam f(D_i)^\delta=\sum_{i\in I_g(\alpha,\delta)}\diam f(D_i)^\delta+\sum_{i\in I_b(\alpha,\delta)}\diam f(D_i)^\delta$$
As soon as $\alpha<\delta<1$,  estimate (\ref{football}) gives 
$$\sum_{i\in I_g(\alpha,\delta)}\diam f(D_i)^\delta\leq C\sum_{i\in I_g(\alpha,\delta)}\diam(D_i)^\alpha$$
Furthermore, by Lemma \ref{hardpart}, there exists an exponent $ \alpha < \delta < 1-\gamma_1 \,\varepsilon^2$ such that the series
$$\sum_{i\in I_b(\alpha,\delta)}\diam f(D_i)^\delta$$
is bounded independently of the covering $D_j$. Thus the entire sum $\sum_i\diam(f(D_i))^\delta$
remains bounded as $\sup_i\diam(D_i)\to 0$. This means  $\dim f(A)\leq\delta\leq 1-\gamma_1\,\varepsilon^2$, and completes  the proof of Theorem \ref{mainlemma}.
 \end{proof} 
 \medskip

By symmetry, c.f. (\ref{factors}),  Theorem \ref{mainlemma}  proves bounds also for expansion of dimension.

\begin{coro}\label{maincoro}
There are constants $c_0, \, c_1 > 0$ such that if $E\subset\partial\D$ and $f:\C\rightarrow\C$ is $K$-quasiconformal with $K=1+\varepsilon$, then
$$\dim(E)\leq 1-c_0\,\varepsilon^2\,\,\,\Rightarrow\,\,\,\dim(f(E))<1-c_1\,\varepsilon^2$$
when $\varepsilon>0$ is small enough.
\end{coro}

Very recently, I. Prause \cite{Pr} has obtained  a different proof for  Theorem \ref{ensi} and Corollary \ref{maincoro}, based in the ideas on \cite{A} and a well known result from Becker and Pommerenke \cite{BP}  which says that
\begin{equation} \label{pomm}
\dim(\Gamma)\leq 1+37\left(\frac{K-1}{K+1}\right)^2
\end{equation}
for every $K$-quasicircle $\Gamma$.  

\begin{rem} \label{yla} Similarly as the  compression bound (\ref{squareddistortion})  led to Theorem \ref{maintheorem}, the inequality (\ref{pomm}) yields  improved upper estimates.
We have hence the symmetric strict inequalities: 

If  $\phi:\C\rightarrow\C$ is a $K$-quasiconformal mapping and $E\subset\partial\D$  with  $\dim(E)=1$, then 
$$\frac{2}{K+1} < \dim(\phi(E)) < \frac{2K}{K+1}.$$
\end{rem}
\bigskip

\noindent Moreover, for the dimension of quasicircles  Smirnov  (unpublished) has obtained the upper bound
$$\dim(\Gamma)\leq 1+\left(\frac{K-1}{K+1}\right)^2,$$
answering a question in \cite{A}. It is still unknown if this bound is sharp; the best known lower bounds so far \cite{ARS} give curves with dimension $1+0.69\left(\frac{K-1}{K+1}\right)^2$. 

The arguments we have used  are related to the {\it{generalized Brennan conjecture}}, which says that
\begin{equation}\label{genb}
\beta_f(p)\leq \frac{p^2}{4}\,\left(\frac{K-1}{K+1}\right)^2 \quad \quad {\rm for } \;\; |p| \leq 2 \frac{K+1}{K-1},
\end{equation}
whenever $f$ is conformal in $\D$ and admits a $K$-quasiconformal extension to $\C$.
This connection suggests the following

\begin{quest}\label{q}
Let $E\subset\R$ be a set with Hausdorff dimension $1$, and let $\phi$ be a $K$-quasiconformal mapping. Is it true that then
\begin{equation}
1-\left(\frac{K-1}{K+1}\right)^2\leq\dim(\phi(E))\leq 1+\left(\frac{K-1}{K+1}\right)^2
\end{equation}
\end{quest}

The positive answer for the right hand side inequality follows from Smirnov's unpublished work, while the left hand side is only known  up to some multiplicative constants. 
On the other hand, Prause \cite{Pr}  proves  the left inequality for the mappings that preserve the unit circle $\partial \D$.

\section{Improved Painlev\'e Theorems}

A compact set $E$ is said to be {\it{removable for bounded analytic functions}} if for any open set $\Omega$ with $E\subset\Omega$, every bounded analytic function on $\Omega\setminus E$ has an analytic extension to $\Omega$. Equivalently, such sets are described by the condition $\gamma(E)=0$, where  $\gamma$ is the analytic capacity 
$$\gamma(E)=\sup\{|f'(\infty)|: f\in H^\infty(\C\setminus E), f(\infty)=0, \|f\|_\infty=1\}$$
Finding a geometric characterization for the sets of zero analytic capacity was a long standing problem.  It was  solved by G. David \cite{D} for sets of finite length, and finally by X. Tolsa \cite{T1} in the general case. The difficulties of dealing with this question motivated the study of  related problems. In particular, we have the question of determining the {\it{removable sets for BMO analytic functions}}, that is, those compact sets $E$ such that every $BMO$ function in the plane, holomorphic on $\C\setminus E$, admits an entire extension. This problem was solved by Kaufman (see \cite{K}), who showed that a set $E$ has this BMO-removability property if and only if $\H^1(E)=0$. 

For the original case of bounded functions the  Painlev\'e condition $\H^1(E)=0$ can be weakened. As is well known, there are sets $E$ with zero analytic capacity and positive length (see \cite{G} for an example). In fact, it is now known that among the compact sets $E$ with $0<\H^1(E)<\infty$, precisely the purely unrectifiable ones are the removable sets for bounded analytic functions \cite{D}. Moreover, if $E$ has positive $\sigma$-finite length, this characterization still remains  true, due to the countable semiadditivity of analytic capacity \cite{T1}.\\
\\
The preceeding problems can  be formulated also in the $K$-quasiregular setting. More precisely, a set $E$ is said to be {\it{removable for bounded (resp. $BMO$) $K$-quasiregular mappings}}, 
if every $K$-quasiregular mapping in $\C\setminus E$ which is in $L^\infty(\C)$ $\left({\rm resp.} BMO(\C)\,  \right)$ admits a $K$-quasiregular extension to $\C$. For simplicity, we use here the term {\it K-removable}  for sets that are removable for bounded $K$-quasiregular mappings.\\
\\
Obviously, when $K=1$, in both situations  $L^\infty$ and $BMO$ we recover the original analytic problem. Moreover, by means of the Stoilow factorization, one can represent any bounded $K$-quasiregular function as a composition of  a bounded analytic function and a  $K$-quasiconformal mapping. The corresponding result holds true also for $BMO$ since  this space, like  $L^\infty$, is  quasiconformally invariant. 

Therefore, when we ask ourselves if a set $E$ is $K$-removable, we just need  to analyze how it may be distorted under quasiconformal mappings, and then apply the  known results for the analytic situation. With this basic scheme, it is shown in \cite[Corollary 1.5]{A}  that every set with dimension strictly below $\frac{2}{K+1}$ is $K$-removable. Indeed, the precise  formulas for the distortion of dimension (\ref{distortionofdimension}) ensure that for such sets the $K$-quasiconformal images have dimension strictly smaller than $1$. 

 Iwaniec and Martin \cite{IM} had earlier conjectured that, more generally,  sets of zero $\frac{2}{K+1}$-dimensional measure are $K$-removable. A preliminary answer to this question was found in  \cite{AIM}, and actually  it was that argument which suggested Theorem \ref{zerotozero}. Using our results from above we can now prove that  sets of zero $\frac{2}{K+1}$-dimensional  measure are even $BMO$-removable.

\begin{coro}\label{bmoremo}
Let $E$ be a compact subset of the plane. Assume that $\H^\frac{2}{K+1}(E)=0$. Then $E$ is removable for all $BMO$ $K$-quasiregular mappings.
\end{coro}
\begin{proof}
Assume that $f\in BMO(\C)$ is $K$-quasiregular on $\C\setminus E$. Denote by $\mu$ the Beltrami coefficient of $f$, and let $\phi$ be the principal solution to $\overline\partial\phi=\mu\,\partial\phi$. Then, $F=f\circ\phi^{-1}$ is holomorphic on $\C\setminus\phi(E)$ and $F\in BMO(\C)$. On the other hand, as we showed in Theorem  \ref{zerotozero}, $\H^1(\phi(E))=0$. Thus, $\phi(E)$ is a removable set for $BMO$ analytic functions. In particular, $F$ admits an entire extension and  $f=F\circ\phi$ extends quasiregularly to the whole plane. 
\end{proof}

We believe that Corollary \ref{bmoremo} is sharp, in the sense that we expect a positive answer to the following 

\begin{quest}
Does there exist for every $K \geq 1$ a compact  set $E$ with $0 < \H^\frac{2}{K+1}(E) < \infty$, such that $E\,$  is not removable for  some $K$-quasiregular functions in $BMO(\C)$.
\end{quest}

Here we observe that by \cite[Corollary 1.5]{A},  for every $t>\frac{2}{K+1}$ there exists a compact set $E$ with dimension $t$, nonremovable for bounded and hence in particular nonremovable for $BMO$ $K$-quasiregular mappings. \\
\\
Next  we return  back to the problem of removable sets for bounded $K$-quasiregular mappings. Here Theorem \ref{zerotozero} proves the conjecture of Iwaniec and Martin that sets with $\H^\frac{2}{K+1}(E) =0$ are $K$-removable. However, the  analytic capacity is somewhat smaller than length, and hence with  Theorem \ref{finitetofinite} we may go even further: If a set has finite or $\sigma$-finite $\frac{2}{K+1}$-measure, then all $K$-quasiconformal images of $E$ have at most $\sigma$-finite length. Such images  may still be removable for  bounded analytic  functions, if we can make sure that the rectifiable part of these sets has zero length. But for this Theorem \ref{maintheorem} provides exactly  the correct tools. We end up with the following improved version of Painlev\'es theorem for quasiregular mappings.

\begin{theo}
Let $E$ be a compact set in the plane, and let $K>1$. Assume that $\H^\frac{2}{K+1}(E)$ is $\sigma$-finite. Then $E$ is removable for all bounded $K$-quasiregular mappings.

In particular, for any $K$-quasiconformal mapping $\phi$ the image $\phi(E)$ is purely unrectifiable.
\end{theo}
\begin{proof}
Let $f:\C\rightarrow\C$ be bounded, and assume that $f$ is $K$-quasiregular on $\C\setminus E$. As in
Corollary \ref{bmoremo} we may find  the principal quasiconformal homeomorphism $\phi: \C \to \C$, such that $F = f \circ \phi^{-1}$ is analytic in $\C \setminus \phi(E)$.  If we can extend $F$ holomorphically to the whole plane,  we are done. Thus  we have to show that $\phi(E)$ has zero analytic capacity.\\
By Theorem \ref{finitetofinite}, $\phi(E)$ has $\sigma$-finite length, that is, $\phi(E)=\cup_nF_n$ where each $\H^1(F_n) < \infty$. A well known result due to Besicovitch (see e.g.\cite[p.205]{M}) assures that each set $F_n$ can be decomposed as 
$$F_n=R_n\cup U_n\cup B_n$$
where $R_n$ is a $1$-rectifiable set, $U_n$ is a purely $1$-unrectifiable set, and $B_n$ is a set of zero length. Because of the semiadditivity of analytic capacity \cite{T1},
$$\gamma(F_n)\leq C\,(\gamma(R_n)+\gamma(U_n)+\gamma(B_n))$$
Now, $\gamma(B_n)\leq  C\, \H^1(B_n)=0$ and $\gamma(U_n)=0$ since purely $1$-unrectifiable sets of finite length have zero analytic capacity \cite{D}. On the other hand, $R_n$ is a $1$-rectifiable image, under a $K$-quasiconformal mapping, of a set of dimension $\frac{2}{K+1}$. Thus applying Theorem \ref{maintheorem} and Corollary \ref{maintheorem4}  to $\phi^{-1}$ shows that  we must have $\H^1(R_n)=0$. Therefore we get $\gamma(F_n)=0$ for each $n$. Again by countable semiadditivity of analytic capacity we conclude  $\gamma(\phi(E))=0$.
\end{proof}

As pointed out earlier, the above theorem does not hold for $K=1$. Any $1$-rectifiable set such as  $E=[0,1]$ of finite and positive length gives a counterexample.  
 In the above proof the  improved distortion of $1$-rectifiable sets was the decisive phenomenon allowing the result. In fact, such good behavior of rectifiable sets has further consequences. 
 For instance, even  {\it strictly above} the critical dimension $\frac{2}{K+1}=1-\frac{K-1}{K+1}$ one may find removable sets, as soon as they have enough geometric regularity.

\begin{coro}\label{more}
There exists  a constant $c\geq 1$ such that if $E\subset\partial\D$ is compact and
$$\dim(E) < 1-c\,\left(\frac{K-1}{K+1}\right)^2$$
then $E$ is removable for bounded and BMO  $K$-quasiregular mappings, $K=1+\varepsilon$, whenever $\varepsilon>0$ is small enough.
\end{coro}
\begin{proof}
This is a consequence of Corollary \ref{maincoro}. If $\varepsilon>0$ is small enough and $K=1+\varepsilon$, then the $K$-quasiconformal images of $E$ will always have dimension strictly below $1$, so that $\gamma(\phi(E))=0$ for each $K$-quasiconformal mapping $\phi$.
\end{proof}

In conjunction with Question \ref{q} we have 

\begin{quest} Let $K>1$. Is then every set $E \subset \partial \D$ with $\dim(E) < 1 - \left(\frac{K-1}{K+1}\right)^2$ removable for bounded and $BMO$ $K$-quasiregular mappings ±
\end{quest}

\section{Examples of extremal distortion}

The previous sections provide a delicate analysis of distortion of $1$-dimensional sets under quasiconformal mappings but  still leave open  the cases  where $\dim(E)=\frac{2}{K+1}$ precisely but $E\,$ 
does not have $\sigma$-finite $\frac{2}{K+1}$-measure. Hence we are faced with the natural question:  Are  there compact sets $E$, with $\dim(E)=\frac{2}{K+1}$, that are non removable for some bounded $K$-quasiregular mappings.

 In this last section  we give a positive answer and show that our  results are sharp in a quite strong  sense. Indeed, to compare with the analytic removability recall first that by Mattila's theorem  \cite{M2},  if a compact set $E$ supports a probability measure  with
 $\mu(B(z,r)) \leq r\, \varepsilon(r)\,  $ and
 \begin{equation} \label{kesa}
  \int_0 \frac{\varepsilon(t)^2}{t} dt < \infty,
 \end{equation}
 then the analytic capacity $\gamma(E) > 0$. On the other hand, if the integral in (\ref{kesa}) diverges, then there are  compact sets $E$ of vanishing analytic capacity  supporting a probability measure with $\mu(B(z,r)) \leq r\, \varepsilon(r)\,  $ \cite{T1}. In a complete analogy we prove

\begin{theo}\label{construction}
Let $K\geq 1$. Suppose  $h(t)=t^\frac{2}{K+1}\,\varepsilon(t)$ is a measure function such that
 \begin{equation} \label{kesa2}
  \int_0 \frac{\varepsilon(t)^{1+1/K}}{t} dt < \infty
 \end{equation}
 Then there is a compact set $E\, $ which is not $K$-removable and yet supports a probability measure $\mu$, with $\mu(B(z,r)) \leq h(r)$ for every $z$ and $r>0$.
 \smallskip
 
 In particular, whenever  $\varepsilon(t)$ is chosen so that in addition for every  $\alpha >0$ we have  $t^\alpha/\varepsilon(t) \to 0$  as  $t \to 0$,   
then  the construction gives  a non-$K$-removable set $E$ with $\dim(E)=\frac{2}{K+1}$.
\end{theo}
\begin{proof} 
We will construct a compact set $E$ and a $K$-quasiconformal mapping $\phi$ such that $\H^h(E)\simeq 1$, and at the same time $\phi(E)$ has a positive and finite $\H^{h'}$-measure for some measure function $h'(t)=t\,\varepsilon'(t)$ where
$$h'(t)=t\,\varepsilon'(t) \quad \mbox{ with} \quad \int_0^1\frac{\varepsilon'(t)^2}{t}dt < \infty $$
Mattilas theorem  shows then $\gamma(\phi(E))>0$, so that there exists  non-constant bounded functions $h$ holomorphic on $\C\setminus\phi(E)$. Thus with $f=h\circ\phi$ we see that $E$ is not removable for bounded $K$-quasiregular mappings.

We will construct the $K$-quasiconformal mapping $\phi$ as the limit of a sequence $\phi_N$ of $K$-quasiconformal mappings, and $E$ will be a Cantor-type set. To reach the optimal estimates we need to change,  at every step in the construction of $E$, both  the size and  the number  $m_j$ of the generating disks. 

Without loss of generality we may assume  that for every  $\alpha >0$,  $t^\alpha/\varepsilon(t) \to 0$  as  $t \to 0$.
\\
{\bf{Step 1}}. Choose first $m_1$ disjoint disks $D(z_i,R_1) \subset \D$, $i=1,...,m_1$, 
so that 
$$c_1:=m_1\,R_1^2 \in (\frac{1}{2}, 1)$$
For $R_1$ small enough (i.e. for $m_1$ large enough) this is clearly possible. The function $f(t)=m_1\,h(t R_1)$ is continuous with  $f(0)=0$. Moreover, for each fixed $t$
$$f(t)=m_1 (tR_1)^{\frac{2}{K+1}}\,\varepsilon(tR_1)=
\frac{\varepsilon\left(t  \sqrt{c_1 / m_1}\,\right)}{\left(t  \sqrt{c_1 / m_1}\right)^{\frac{2}{K+1}} }\; t^2 c_1 \; \to \;  \infty$$
as $m_1\to\infty $.
Hence for any $t<1$ we may choose $m_1$ so large that there exists $\sigma_1\in(0,t)$ satisfying $m_1\,h(\sigma_1^K R_1)=1$.  A simple calculation gives
\begin{equation}
m_1\; \sigma_1 R_1 \,\varepsilon(\sigma_1^K R_1)^\frac{K+1}{2 K}\,(c_1)^\frac{1-K}{2K}=1
\end{equation}

Next, let  $r_1=R_1$. For each $i=1,\dots, m_1$, let $\varphi_i^1(z)=z_i+\sigma_1^KR_1\,z$ and, using the notation $\alpha D(z,\rho):= D(z,\alpha \rho)$, set
$$\aligned
D_i&:=\frac{1}{\sigma_1^K}\,\varphi^1_i(\D)=D(z_i, r_1)\\
D'_i&:=\varphi^1_i(\D)=D(z_i,\sigma_1^K r_1) \subset D_i
\endaligned$$ 
As  the first approximation of the mapping define
$$
 g_1(z)=
\begin{cases}
\sigma_1^{1-K}(z-z_i)+z_i, &z\in D'_i\\
\left|\frac{z-z_i}{r_1}\right|^{\frac{1}{K}-1}(z-z_i)+z_i, \; &z\in D_i\setminus D'_i\\
z, & z \notin \cup D_i
\end{cases}
$$
This is a $K$-quasiconformal mapping, conformal outside of $\bigcup_{i=1}^{m_1}(D_i\setminus D'_i)$. It maps each $D_i$ onto itself and $D'_i$ onto $D''_i=D(z_i,\sigma_1r_1)$, while the rest of the plane remains fixed. Write $\phi_1=g_1$. 
\\
{\bf{Step 2}}. We have already fixed $m_1, R_1, \sigma_1$ and $c_1$. Consider $m_2$ disjoint disks of radius $R_2$, centered at $z^2_j$, $j=1,\dots, m_2$, uniformly distributed inside of $\D$, so that
$$c_2=m_2\,R_2^2>\frac{1}{2}$$
Then repeat the above procedure and choose $m_2$ so large that the equation 
$$m_1\, m_2\;  h(\sigma_1^K \sigma_2 ^K R_1R_2)=1$$
has a unique solution $\sigma_2\in(0,1)$, as small as we wish. Then,
$$m_1m_2\,\sigma_1 \sigma_2 R_1R_2 \,\varepsilon(\sigma_1^K \sigma_2^K R_1R_2)^\frac{K+1}{2K}\,(c_1c_2)^\frac{1-K}{2K}=1$$
Denote $r_2=R_2\,\sigma_1r_1$ and $\varphi^2_j(z)=z^2_j+\sigma_2^KR_2\,z$, and define the auxiliary disks
$$\aligned
D_{ij}=\phi_1\left(\frac{1}{\sigma_2^K}\, \varphi^1_i \circ \varphi^2_j(\D)\right)=D(z_{ij}, r_2)\\
D'_{ij}=\phi_1\left(\, \varphi^1_i \circ \varphi^2_j(\D)\right)=D'(z_{ij}, \sigma_2^K r_2)
\endaligned$$
for certain $z_{ij} \in \D$, where $i=1,\dots,m_1$ and $j=1,\dots,m_2$. Now Let
$$
g_2(z)=
\begin{cases}
\sigma_2^{1-K}(z-z_{ij})+z_{ij}&z\in D'_{ij}\\
\left|\frac{z-z_{ij}}{r_2}\right|^{\frac{1}{K}-1}(z-z_{ij})+z_{ij}&z\in D_{ij}\setminus D'_{ij}\\
z&\text{otherwise}
\end{cases}
$$
Clearly, $g_2$ is $K$-quasiconformal, conformal outside of $\bigcup_{i,j}(D_{ij}\setminus D'_{ij})$, maps each $D_{ij}$ onto itself and $D'_{ij}$ onto $D''_{ij}=D(z_{ij},\sigma_2r_2)$, while the rest of the plane remains fixed. 
Define $\phi_2=g_2\circ\phi_1$.\\

\begin{figure}[ht]
\begin{center}
\includegraphics{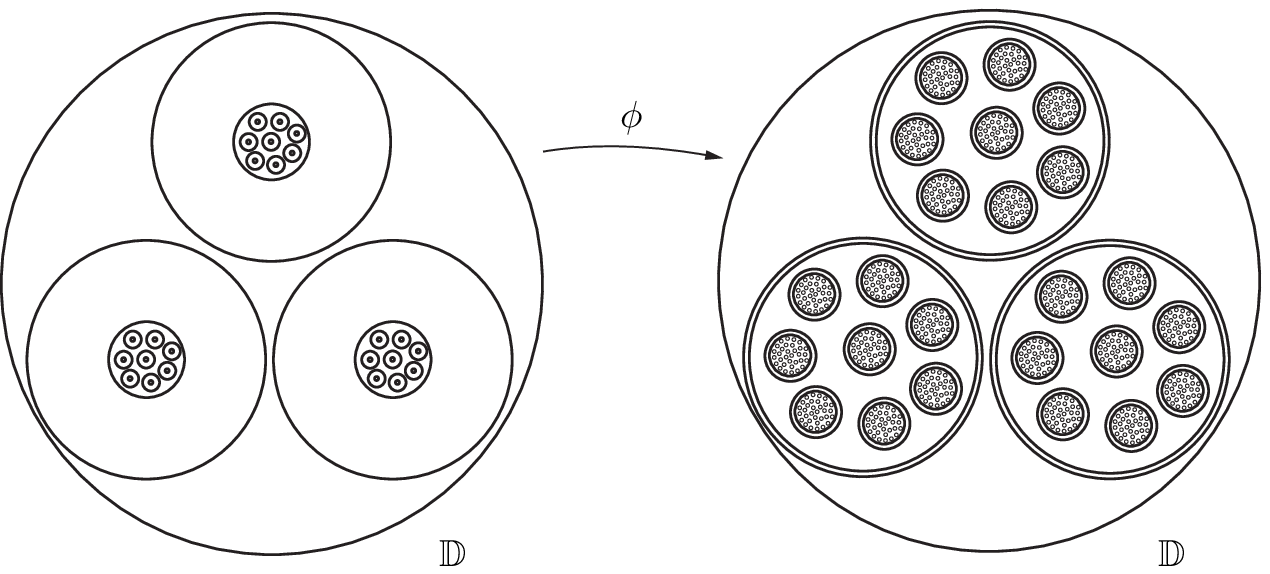}
\end{center}
\end{figure}

\noindent  {\bf{The induction step}}. After step $N-1$ we take $m_N$  disjoint disks of radius $R_N$, with union of $D(z^N_l, R_N)$ covering  at least half of the area of $\D$,
\begin{equation}\label{cc}
c_N=m_N\,R_N^2>\frac{1}{2}
\end{equation}
As before  we may choose $m_N$ so large  that 
$m_1\dots m_N\,h(\sigma_1^K \dots\sigma_N^K \,R_1\dots R_N)=1$
holds for a unique  $\sigma_N$, as small as we wish. Note that $\displaystyle\lim_{N\to\infty}\sigma_N=0$ and
$$m_1\dots m_N\, \sigma_1 R_1\dots \sigma_N R_N \,\varepsilon(\sigma_1^K R_1\dots \sigma_ N^K R_N)^\frac{K+1}{2K}(c_1\dots c_N)^\frac{1-K}{2K}=1$$
Denote then $\varphi^N_j(z)=z^N_j+\sigma_N^KR_N\,z$ and $r_N=R_N\,\sigma_{N-1}r_{N-1}$. For any multiindex $J=(j_1,...,j_N)$, where $1\leq j_k\leq m_k$, $k=1,...,N$, let 
$$\aligned
D_{J}=\phi_{N-1}\left(\frac{1}{\sigma_N^K}\, \varphi^1_{j_1}\circ\dots\circ\varphi^N_{j_N}(\D)\right)=D(z_{J}, r_N)\\
D'_{J}=\phi_{N-1}\left(\varphi^1_{j_1}\circ\dots\circ\varphi^N_{j_N}(\D)\right)=D'(z_{J}, \sigma_N^K r_N)
\endaligned$$
and let
$$
g_N(z)=
\begin{cases}
\sigma_N^{1-K}(z-z_{J})+z_{J}&z\in D'_{J}\\
\left|\frac{z-z_{J}}{r_N}\right|^{\frac{1}{K}-1}(z-z_{J})+z_{J}&z\in D_{J}\setminus D'_{J}\\
z&\text{otherwise}
\end{cases}
$$
Clearly, $g_N$ is $K$-quasiconformal, conformal outside of $\bigcup_{J=(j_1,...,j_N)}(D_{J}\setminus D'_{J})$, maps $D_{J}$ onto itself and $D'_{J}$ onto $D''_{J}=D(z_{J},\sigma_Nr_N)$, while the rest of the plane remains fixed. 
Now define $\phi_N=g_N\circ\phi_{N-1}$.\\  
\\
Since each $\phi_N$ is $K$-quasiconformal and  equals the identity  outside the unit disk $\D$, there exists a limit $K$-quasiconformal mapping  
$$\phi=\lim_{N\to\infty}\phi_N$$
with convergence in $W^{1,p}_{loc}(\C)$ for any $p<\frac{2K}{K-1}$. On the other hand, $\phi$ maps the compact set $$E=\displaystyle\bigcap_{N=1}^\infty\left(\bigcup_{j_1,...,j_N}\varphi^1_{j_1}\circ\dots\circ\varphi^N_{j_N}(\D)\right)$$
to the set
$$\phi(E)=\bigcap_{N=1}^\infty\left(\bigcup_{j_1,...,j_N}\psi^1_{j_1}\circ\dots\circ\psi^N_{j_N}(\D)\right)$$
where  we have written $\psi^i_j(z)=z^i_j+\sigma_iR_i z$, $j=1,...,m_i$, $i \in \N$. 

To complete the proof, write
\begin{equation}\label{radii}
s_N=(\sigma_1^KR_1)\dots(\sigma_N^KR_N) \quad \mbox{ and } \quad
t_N=(\sigma_1R_1)\dots(\sigma_NR_N)
\end{equation}
Observe that  we have chosen the parameters $R_N, m_N, \sigma_N$  so that
\begin{eqnarray}
&&m_1\dots m_N\,h(s_N)=1\label{viides}\\
&&m_1\dots m_N\,t_N\,\varepsilon(s_N)^\frac{K+1}{2K}(c_1\dots c_N)^\frac{1-K}{2 K}=1 \label{toinen}
\end{eqnarray}

{\bf{Claim}}. $\H^{h}(E)\simeq 1$. \\
Since 
 $\diam(\varphi^{1}_{j_{1}}\circ\dots\circ\varphi^N_{j_N}(\D))\leq \delta_N \to 0$ when $N\to \infty$, we have by (\ref{viides}) 
$$\H^h(E)=\lim_{\delta\to 0}\H^{h}_{\delta}(E)\leq \lim_{\delta\to 0}\sum_{j_1,\dots,j_{N}}h(\diam(\varphi^{1}_{j_{1}}\circ\dots\circ\varphi^N_{j_N}(\D)))=m_1\dots m_N\,h(s_N)=1$$
For the converse inequality, take a finite covering $(U_j)$ of $E$ by open disks of diameter $\diam(U_j)\leq\delta$ and let $\delta_0=\inf_j(\diam(U_j))>0$.
Denote by $N_0$ the minimal integer such that $s_{N_0}\leq\delta_0$. By construction, the family $(\varphi^{N_0}_{j_{N_0}}\circ\dots\circ\varphi^1_{j_1}(\D))_{j_1,\dots,j_{N_0}}$ is a covering of $E$ with the $\M^{h}$-packing condition \cite{M}. Thus, 
$$\aligned
\sum_jh(\diam(U_j))&\geq C\,\sum_{j_1,\dots,j_{N_0}}h(\diam(\varphi^{N_0}_{j_{N_0}}\circ\dots\circ\varphi^1_{j_1}(\D)))=C\endaligned$$
Hence, $\H^{h}_{\delta}(E)\geq C$ and 
 letting $\delta\to 0$, we get that
$$C\leq\H^{h}(\phi(E))\leq 1$$
proving our first claim. \\

A similar argument, based this time on (\ref{toinen}), gives that $\H^{h'}(\phi(E))\simeq 1$ for a  measure function $h'(t) = t \varepsilon'(t)$, as soon as for all indexes $N$
\begin{equation} \label{kolmas}
 \varepsilon'(t_N) = \varepsilon(s_N)^\frac{K+1}{2K}(c_1\dots c_N)^\frac{1-K}{2 K}
\end{equation}

{\bf{Claim}}. One can find a continuous and nondecreasing function $\varepsilon'(t)$ satisfying (\ref{kolmas}) and
\begin{equation} \label{neljäs}
\int_0^1\frac{\varepsilon'(t)^2}{t}dt < \infty
\end{equation}
Indeed, let us first choose a continuous nondecreasing function $v(t)$ so that $v(t)\to 0$ as $t\to0$ and  so that (\ref{kesa2})  still holds in the  form
\begin{equation}
  \int_0 \frac{\varepsilon(t)^{1+1/K}}{ t \, v(t) } dt < \infty
\end{equation}
 In the above inductive construction we can then  choose  the  $\sigma_j$'s so that
 $v(\sigma_1^K \cdots \sigma_N^K) \leq 2^{-N (1-1/K)}$ for every index $N$. Now  (\ref{cc}) and  (\ref{kolmas})  imply $$\varepsilon'(t_N)^2 \leq \varepsilon(s_N)^{1+1/K} \, 2^{N (1-1/K)} \leq \frac{\varepsilon(s_N)^{1+1/K}}{v(s_N)}$$
On the other hand by (\ref{radii}) we also have $t_{N-1}/t_N \leq s_{N-1}/s_N  $ and so we may extend $\varepsilon'(t)$, determined by (\ref{kolmas})  only at the $t_N$'s, so that it is continuous, nondecreasing and satisfies
$$  \int_0 \varepsilon'(t)^{2} \; \frac{dt}{ t  } \leq   \int_0 \frac{\varepsilon(s)^{1+1/K}}{  v(s) }\; \frac{ ds}{s} < \infty
$$
Hence the claim follows. Combining it with Mattila's theorem \cite{M2}  completes the proof of the Theorem.
\end{proof}
\bigskip

Lastly let us note that if we do not care for the analytic capacity of the target set, a straightforward modification of the previous Theorem, normalizing the disks of the construction so that
$m_N\,  t_N\,  \eta(t_N) = 1$, gives 
\begin{coro} Let $K \geq 1$ and let $h(t) = t \, \eta(t)$ be a measure function such that  
\begin{itemize}
\item $\eta$ is continuous and nondecreasing, $\eta(0)=0$ and $\eta(t)=1$ whenever $t\geq 1$.
\item $\displaystyle\lim_{t\to 0}\frac{t^\alpha}{\eta(t)}=0$ for all $\alpha>0$.
\end{itemize}
There exists a compact set $E\subset\D$ and a $K$-quasiconformal mapping $\phi$ such that
\begin{equation}
\dim(E)=\frac{2}{K+1} \quad {\rm and } \quad \H^h(\phi(E))=1
\end{equation}
\end{coro}

\vskip 1cm
\begin{itemize}
\item[]{Department of Mathematics and Statistics, University of Helsinki, P.O.
Box 68, FIN-00014, University of Helsinki, Finland\\
{\it E-mail address:} { kari.astala@helsinki.fi}}

\item[]{Departament de Matem\`atiques, Facultat de Ci\`encies, Universitat Aut\`onoma de Barcelona, 08193-Bellaterra, Barcelona, Catalonia\\
{\it E-mail address:} {albertcp@mat.uab.cat}, {mateu@mat.uab.cat}, {orobitg@mat.uab.cat}}

\item[]{Mathematics Department, 202 Mathematical Sciences Bldg., University of Missouri, Columbia, MO 65211-4100, USA 
{\it E-mail address:} {ignacio@math.missouri.edu}}
\end{itemize}

\end{document}